\documentclass[12pt]{amsart}

\usepackage{amsmath,amsthm,amsfonts,latexsym,euscript}

\newtheorem{df}{Definition}[section]
\newtheorem{lm}[df]{Lemma}
\newtheorem{prop}[df]{Proposition}
\newtheorem{thm}[df]{Theorem}
\newtheorem{cor}[df]{Corollary}
\newtheorem{rmk}[df]{Remark}
\newtheorem{rmks}[df]{Remarks}
\newtheorem{prob}[df]{Problem}

\newcommand{\bgdf}{\begin{df}}
\newcommand{\nddf}{\end{df}}
\newcommand{\bglm}{\begin{lm}}
\newcommand{\ndlm}{\end{lm}}
\newcommand{\bgprop}{\begin{prop}}
\newcommand{\ndprop}{\end{prop}}
\newcommand{\bgthm}{\begin{thm}}
\newcommand{\ndthm}{\end{thm}}
\newcommand{\bgcor}{\begin{cor}}
\newcommand{\ndcor}{\end{cor}}
\newcommand{\bgrmk}{\begin{rmk}}
\newcommand{\ndrmk}{\end{rmk}}
\newcommand{\bgrmks}{\begin{rmks}}
\newcommand{\ndrmks}{\end{rmks}}
\newcommand{\bgprob}{\begin{prop}}
\newcommand{\ndprob}{\end{prop}}
\newcommand{\bgeq}{\begin{eqnarray}}
\newcommand{\ndeq}{\end{eqnarray}}
\newcommand{\bgeqq}{\begin{eqnarray*}}
\newcommand{\ndeqq}{\end{eqnarray*}}

\newcommand{\pf}{{\em Proof. }} %  \pf is already defined in amslatex
\newcommand{\QED}{\hfill Q.E.D.} % \qed (not \QED) is defined in amslatex
\newcommand{\vv}{\vspace{4mm}\\}

\newcommand{\vvvv}{\vspace{1cm}\\}

\numberwithin{equation}{section}

\begin{document}

\title[Quantum $ax + b$ Group as Quantum Automorphism Group]
{Quantum $a x + b$ Group as Quantum Automorphism Group of $k[x]$
%The Algebraic Level
}
\author[Shuzhou Wang]
{\bf Shuzhou Wang}
\address{Department of Mathematics, University of Georgia,
Athens, GA 30602
\newline \indent
Fax: 706-542-2573
}
\email{szwang@math.uga.edu}
%\dedicatory
%\date{Dec, 1999}
\thanks{Partially supported by NSF grant DMS-9627755.}

 % %Top matter begins here
 %\title[...]{...}
 %\author[...]{...}
 %\address{...}
 %\curraddr{...}
 %\email{...}
 %\urladdr{...}
 %\dedicatory{}
 %\date{\today}
 %\thanks{NSF}
 %\translator{...}
 %\keywords{...}
 %\subjclass{...}
 %\begin{abstract}...\end{abstract}
 % %Top matter ends here
 %\maketitle

%\begin{document}

\begin{abstract}
By introducing a result that guarantees a
given bialgebra to be a Hopf algebra under a natural condition,
we show that the quantum automorphism group of
the algebra $k[x]$ of polynomials %in one variable
over a field $k$ (of any characteristic) is
the universal quantum $a x + b $ group ${\mathcal A}$, generalizing
the fact that the automorphism group of $k[x]$
is the $a x + b$ group. The $q$-deformation of the
$a x + b$ group is then seen as one among 
a certain family ${\mathcal A}_{q, n}$ ($q \in k$, $n$ an integer) 
of quantum subgroups of this universal quantum group. % ${\mathcal  A}$.
%We also discuss other quantum subgroups of ${\mathcal  A}$.
\end{abstract}

\maketitle

%\vspace{8mm}
%\noindent
\section{Introduction}

In \cite{W15}, adapting the philosophy that quantum groups should
be viewed as the mathematical objects encoding quantum symmetries
of (noncommutative) spaces \cite{Mn1}, instead of as deformations of
ordinary groups, we described the quantum automorphism groups
of matrix algebras over $\Bbb C$. Just as automorphism
groups of matrix algebras over $\Bbb C$ are compact Lie groups,
the quantum automorphism groups of such algebras
are compact matrix quantum groups in the sense of Woronowicz
(see \cite{Wor5,W1,W5,W5',W7,W15}). These quantum automorphism groups
are of a quite different nature from the $q$-deformations of Lie groups, 
although their representation ring is the same as the representation ring 
of $SO(3)$ for the generic situation, 
as cleverly shown by Banica \cite{Banica7}. 
These quantum automorphism groups contain 
the ordinary automorphism (Lie) groups as
{\em proper} subgroups. Hence, quantum group symmetries 
are significantly richer than ordinary group symmetries.

In this note, adapting the same philosophy as in \cite{Mn1,W15} and using  
``algebra of functions on quantum group approach'',
we determine explicitly the quantum automorphism
group of the polynomial algebra $k[x]$ of one variable over a field $k$
(of any characteristic).
It is well known that the automorphism group of the algebra
$k[x]$ is the $a x + b$ group over the field $k$. Namely, 
every automorphism $\alpha$ of the $k$-algebra $k[x]$ is of the form 
$$\alpha(x) = ax + b$$
 for some $a, \; b \in k$ with $a \neq 0$. 
In exactly the same way, we find that the quantum automorphism group of
$k[x]$ is the universal quantum $a x + b$ group, which has been
studied earlier by Sweedler \cite{Sweedler} in a different context.
As a matrix quantum group, its algebra ${\mathcal  A}$ of functions has
generators $a, a^{-1}, b$ ($a$ and $b$ are free with respect to each other)
and is represented by
\bgeqq   T = %(t_{ij})
\left(
\begin{array}{cc}
a & 0    \\
b & 1
\end{array}
\right).
\ndeqq
We only describe this quantum group at the purely algebraic
level in this note, its $C^*$-algebraic (topological) description
(e.g. when $k = {\Bbb R}$ or ${\Bbb C}$) would require a
rather involved analysis of unbounded operators \cite{Wor8}.
This problem seems to have not been solved yet, although there have been
serious attempts at it.
Note that in \cite{Majid1}, Majid describes a
very big quantum semigroup (bialgebra) for $k[x]$. This quantum semigroup is
not a quantum group in our sense (i.e. it has no antipode),
and it not homogeneous in the sense that it does not preserve
the degrees of the polynomials in $k[x]$.
%\vvvv

\section{From Bialgebras to Hopf Algebras: A Sufficient Condition}
%\vv

In practice, there are a profusion of bialgebras (i.e. quantum
semigroups), cf. \cite{Mn1,W15,Majid1}. 
However, as Drinfeld pointed out in his paper \cite{Dr1}, 
before the discovery of quantum groups there were very few 
non-trivial examples of noncommutative and noncocommutative Hopf 
algebras and it was very hard to find such. 
This is one of the reasons why
the Drinfeld-Jimbo construction of such examples has attracted
much attention in both physics and mathematics. In this section,
we give a natural condition that guarantees a given bialgebra to
be a Hopf algebra. It is an analog of the result of Woronowicz in
\cite{Wor9}. This result is very useful in constructing examples
of Hopf algebras. For instance, the quantum groups in \cite{FRT1}
can be viewed as such. See also \cite{W1,W5,W15}.

For a given algebra ${\mathcal A}$ over a field $k$
(of any characteristic), denote by $M_n({\mathcal A})$ the algebra of
$n \times n$ matrices over $\mathcal A$.
\vv
{\bf Theorem 1.}
{\em
Let $({\mathcal A}, \Delta, \epsilon)$ be a unital bialgebra over a field $k$
that is generated by a multiplicative matrix $u=(u_{ij})_{i,j=1}^n$ 
and relations $R$. That is, ${\mathcal A}$ is 
generated by $u_{ij}$ ($i, j = 1 \cdots n$) and 
relations $R$ such that
$$\Delta(u_{ij}) = \sum_{k=1}^n u_{ik} \otimes u_{kj},
\hspace{1cm} \epsilon(u_{ij}) = \delta_{ij}. $$
Then ${\mathcal A}$ is a Hopf algebra if and only if
$u$ is invertible in $M_n({\mathcal A})$ and the entries of $u^{-1}$ 
satisfy the opposite relations $R^{op}$.   
When this condition is satisfied, the antipode $S$ is given on the
generators $u_{ij}$ by $$S(u) = u^{-1}.$$}
\pf (cf Manin \cite{Mn1}.) 
The necessary condition follows from the antipodal property.

We prove that the condition is sufficient. 
Since the entries of $u^{-1}$ satisfy the opposite relations $R^{op}$ 
we can define a anti-homomorphism $S$ on $\mathcal A$ by the formula 
$$
(S(u_{ij}))_{i,j=1}^n = u^{-1} .
$$
We show that $S$ satisfies the antipodal property 
$$
m (S \otimes id ) \Delta(a) = \epsilon(a) 1 
= 
m (id \otimes S ) \Delta(a) 
$$
for all $a \in {\mathcal A}$. 
From the definition of $S$, 
this is clearly true for $a = u_{ij}$. Since 
$\mathcal A$ is generated by the $u_{ij}$'s 
as an algebra, to prove the antipodal property, we 
show that if $a, b \in {\mathcal A}$ satisfy the antipodal property, 
then so does $ab$. 
Using Sweedler's notation, let 
$$
\Delta(a) = \sum a_{(1)} \otimes a_{(2)}, \hspace{1cm}  
\Delta(b) = \sum b_{(1)} \otimes b_{(2)}. 
$$
Then 
\bgeqq
& & m (S \otimes id ) \Delta(ab) = 
    m (S \otimes id ) \Delta(a) \Delta(b) \\
&=& m (S \otimes id ) ( \sum a_{(1)} b_{(1)} \otimes a_{(2)} b_{(2)} ) \\
&=& m \sum S(b_{(1)}) S(a_{(1)}) \otimes a_{(2)} b_{(2)} \\
&=&   \sum S(b_{(1)}) S(a_{(1)}) a_{(2)} b_{(2)} 
=   \sum S(b_{(1)}) [S(a_{(1)}) a_{(2)}] b_{(2)} \\
&=&   \sum S(b_{(1)}) \epsilon(a) b_{(2)} 
=   \epsilon(a) \epsilon(b) 1 = 
\epsilon(a b) 1. 
\ndeqq 
Similarly, we have 
$$ m (id  \otimes S ) \Delta(ab) = \epsilon(a b) 1. 
$$
\QED
\vv
{\em Remark 1.}
Let $u$ be a matrix corepresentation of a Hopf algebra
(or a multiplicative matrix in the sense of Manin \cite{Mn1}).
Let $S$ be the antipode of $H$. Then both $u$ and $u^t$ are invertible,
and we have $u^{-1} = S(u)$, but the inverse $(u^t)^{-1}$ is not
equal to $(u^{-1})^t$ in general.
That is, we have in general
$$S(u^t) = S(u)^t = (u^{-1})^t \neq (u^t)^{-1}. $$
%but $S(u^t) \neq (u^t)^{-1}$.
The reader may convince himself of
this by taking the Hopf algebra of $SU_q(2)$ for $q \in (0, 1)$,  or
the universal quantum group $A_u(Q)$ for $Q>0$ but $Q \neq cI_n$
\cite{W5}.

If we restrict to compact quantum groups in the sense of
Woronowicz \cite{Wor5}, then $S(u^t) = (u^t)^{-1}$
if and only if the quantum group is of Kac type, i.e.
if and only if $S^2=id$.  % the Haar state is a trace.
To see this, it suffices to use the fact that $S$ preserves
the involution in this case.
%\vvvv

\section{The Universal Quantum $a x + b$ Group as a Quantum
Automorphism Group}
%\vv

In this section, we describe the quantum automorphism group of
$k[x]$. It turns out that its Hopf algebra is exactly
the same as the one studied earlier by Sweedler in
a different context (p89 of \cite{Sweedler}).

Note we are taking the ``functions-on-group'' point of view. This means
that a {\bf quantum group} is a Hopf algebra (over a field $k$ of any
characteristic) of ``noncommutative functions''.
An {\bf action} of a quantum group on an algebra is defined to be a
right coaction of the Hopf algebra on the algebra.
A quantum group endowed with an action on an algebra
is called a {\bf quantum transformation group} over the algebra.

Let $({\mathcal  A}_1, \alpha_1)$ and
    $({\mathcal  A}_2, \alpha_2)$ be two quantum transformation groups over
    an algebra ${\mathcal  B}$. A {\bf morphism} $\pi$ from
    $({\mathcal  A}_1, \alpha_1)$ to
    $({\mathcal  A}_2, \alpha_2)$ is defined to be a morphism $\pi$ of
Hopf algebras from ${\mathcal  A}_2$ to ${\mathcal  A}_1$ with the property
   $$\alpha_1 = (id \otimes \pi) \alpha_2.$$
We refer the reader to  \cite{W15} for more discussions on
categories of quantum transformation groups.

Let $k[x]_n \subset k[x]$ be the subspace of
polynomials of degrees less or equal to $n$ ($n \geq 0$).
Note that the automorphism group of $k[x]$ leaves
the subspaces $k[x]_n$ invariant. Therefore, it is natural to expect that
its quantum automorphism group also has this property. On the other hand,
it is clear that the class of quantum groups acting on $k[x]$ that
leave $k[x]_n$ invariant is a category (compare with Definition 2.1
of \cite{W15}).
%We call an action of a quantum group homogeneous if $k_n[x]$ is invariant
%under the action.
Let ${\mathcal  C}$ denote this category.
We now show that the category ${\mathcal  C}$ has a universal
object (then it is unique up to isomorphism by abstract nonsense). 
We explicitly describe it
in terms of generators and relations. This universal object is
the quantum automorphism group of $k[x]$ (compare with \cite{W15}).

Let $({\mathcal  A}_0, \alpha_0)$ be an object in ${\mathcal  C}$, 
where $\alpha_0$
denotes the action of the quantum group ${\mathcal  A}_0$ on $k[x]$. Since
$k[x]_1$ is invariant under $\alpha_0$, there are elements
$a_0, b_0 \in {\mathcal  A}_0$ such that
\bgeqq
\alpha_0 (x)  = x \otimes a_0 + 1 \otimes b_0.
\ndeqq
Since (see Section 2 of \cite{Mn1} and Definition 2.1 of \cite{W15})
$$( 1 \otimes \Delta_0 ) \alpha_0 = ( \alpha_0 \otimes 1 ) \alpha_0,
\; \; \;
 (id \otimes \epsilon_0) \alpha_0 = id,$$
where $\Delta_0$ and $\epsilon_0$ are respectively the coproduct
and counit on ${\mathcal  A}_0$, we see that
\bgeqq
\left(
\begin{array}{cc}
a_0 & 0    \\
b_0 & 1
\end{array}
\right)
\ndeqq
is a multiplicative matrix (see Section 2 of \cite{Mn1}). Namely,
letting $T_0 = (t^0_{ij})$ be the above matrix, then
\bgeqq
 \Delta_0 (t^0_{ij}) = \sum_{r=1}^2 t^0_{ir} \otimes  t^0_{rj};
 \; \; \;
 \epsilon_0 (t^0_{ij}) = \delta_{ij}.
\ndeqq
This combined with the antipodal property for ${\mathcal  A}_0$
gives
\bgeqq
\sum_{r =1}^2 S(t^0_{kr}) t^0_{rl} = \delta_{kl},
\; \; \;
\sum_{r =1}^2 t^0_{kr} S(t^0_{rl}) = \delta_{kl}, \hspace{2cm} (*)
\ndeqq
where $S_0$ is the antipode of ${\mathcal  A}_0$, $k, l = 1, 2$.
Simplifying the above relations $(*)$ we obtain
\bgeqq
S_0 (a_0) = a_0^{-1}, \; \; \;
S_0 (b_0) = -  b_0 a_0^{-1}.
\ndeqq

Now we can state  the following
\vv
{\bf Theorem 2.}
{\em
Let ${\mathcal  A}$ be the universal algebra over the field $k$ generated by
$a, a^{-1}, b$, subjecting to the following two relations
\bgeqq
a a^{-1} = 1 = a^{-1} a.
\ndeqq
Then

(1).
${\mathcal  A}$ is a quantum group with Hopf
algebra structure
\bgeqq
\Delta (a) = a \otimes a, \; \;  \;
\Delta (b) = b \otimes a + 1 \otimes b,
\ndeqq
\bgeqq
%S(a) = a^{-1}, \; \;
%S(b) = - b a^{-1}, \; \;
%\epsilon(a) = 1, \; \; \epsilon(b) = 0.
%
 S
\left(
\begin{array}{cc}
a & 0    \\
b & 1
\end{array}
\right)
=
\left(
\begin{array}{cc}
a^{-1}     & 0    \\
- b a^{-1} & 1
\end{array}
\right),
 \; \; \; \;
 \epsilon
\left(
\begin{array}{cc}
a & 0    \\
b & 1
\end{array}
\right)
=
\left(
\begin{array}{cc}
1 & 0    \\
0 & 1
\end{array}
\right).
\ndeqq

(2).
The quantum group ${\mathcal  A}$ admits a natural action $\alpha$ on
the algebra $k[x]$ given by
\bgeqq
\alpha (x) = x \otimes a + 1 \otimes b.
\ndeqq
The pair $({\mathcal  A}, \alpha)$ is the quantum automorphism group of $k[x]$
(i.e. a universal object in the category ${\mathcal  C}$), and it naturally
contains (in the sense of \cite{W1}) the ordinary
automorphism group $Aut(k[x])$ as a subgroup.}
\vv
Note that part (1) of the above theorem  
%The Hopf algebra ${\mathcal  A}$ 
appeared in Sweedler \cite{Sweedler} pp89-90
in a completely different context, where the result 
was stated without proof. 
Although the proof of (1) can be done by straightforwardly verifying the
axioms of a Hopf algebra one by one, it is also illuminating to
prove it by using Theorem 1 above, which we leave to the reader as
an exercise. The proof
of (2) follows easily from the calculations above. Note that the proof
of this theorem is easier than the proofs of the corresponding
results in \cite{W15}, partly because this theorem is of a purely algebraic
nature whereas the results in \cite{W15} concerns analytic aspect as
well (i.e. Woronowicz Hopf $C^*$-algebras).

In view of this theorem, the quantum group ${\mathcal  A}$ can be
called the {\bf universal quantum $ax + b$ group}.  
\vv {\em Remark 2.} The quantum
automorphism group $({\mathcal  A}, \alpha)$ above is the
universal object in the category of all Hopf algebras coacting
from the {\bf right} on the algebra $k[x]$. It is easy to see that
the pair $({\mathcal  A}^{op}, \alpha^{op})$ is the universal
object in the category of all Hopf algebras that coact from the
{\bf left} on the algebra $k[x]$ and leave subspaces $k_n[x]$
invariant, where ${\mathcal  A}^{op}$ is the Hopf algebra that has
the same elements as ${\mathcal  A}$ but with the opposite product
$m^{op}$ and opposite coproduct $\Delta^{op}$ (see \cite{Mn1}),
and $\alpha^{op}$ is defined by \bgeqq \alpha^{op} (x) = a \otimes
x + b \otimes 1. \ndeqq The multiplicative matrix for $({\mathcal
A}^{op}, \alpha^{op})$ is \bgeqq \left(
\begin{array}{cc}
a & b    \\
0 & 1
\end{array}
\right).
\ndeqq
Hence in a sense we obtain the same quantum group with left coaction.
%\vvvv

\section{$q$-Deformation ${\mathcal  A}_{q}$
of the $a x + b$ Group and Other Quantum Subgroups 
${\mathcal  A}_{q, n}$ of ${\mathcal  A}$ }
%\vv

Let $q \in k^*$ be a non-zero element.  
Let ${\mathcal  I}_q$ be the
ideal of ${\mathcal  A}$ generated by $ab - q ba$.
Then one easily verifies that ${\mathcal  I}_q$ is a
Hopf ideal of ${\mathcal  A}$. Let ${\mathcal  A}_q$
be the quotient ${\mathcal  A}/{\mathcal  I}_q$. Then ${\mathcal  A}_q$
is a deformation of the ordinary $ax + b$ group over $k$ in the sense
that ${\mathcal  A}_1$ is the ordinary algebra of coordinate functions
on this group.

More generally, with $q$ as above and $n \in {\mathbb Z}$ an integer, 
let ${\mathcal  I}_{q, n}$ be the
ideal of ${\mathcal  A}$ by $a^n b - q b a^n$.
Then ${\mathcal  I}_{q, n}$ is a
Hopf ideal of ${\mathcal  A}$. To see this, 
$$
\Delta (a^n b - q b a^n) = (a^n b - q b a^n) \otimes a^{n+1} 
+ a^n \otimes (a^n b - q b a^n) , 
$$
$$
\text{hence},  \; \; \; 
\Delta ({\mathcal  I}_{q, n}) 
\subseteq 
{\mathcal  I}_{q, n} \otimes {\mathcal A} + 
 {\mathcal A} \otimes {\mathcal  I}_{q, n} ; 
$$
$$
\epsilon(a^n b - q b a^n) = 0, \; \; \; 
\text{hence}, \; \; \; 
\epsilon( {\mathcal  I}_{q, n}) = 0 ; 
$$
$$
S(a^n b - q b a^n) = a^{-n} (-a^n b + q b a^n ) a^{n-1}  ,
$$
$$
\text{hence}, \; \; \; 
S({\mathcal  I}_{q, n}) \subseteq {\mathcal  I}_{q, n}.  
$$
Let ${\mathcal  A}_{q, n}$
be the quotient ${\mathcal  A}/{\mathcal  I}_{q, n}$. 
Then ${\mathcal  A}_{q, n}$
is a quantum subgroup of the universal quantum group ${\mathcal  A}$,
and of ${\mathcal  A}_{q^m, mn}$ for integers  
$m \neq 0$ (because $a^{m n} b - q^m b a^{m n} = 0$ 
in  ${\mathcal  A}_{q, n}$). 
In particular, since ${\mathcal  A}_{q}= {\mathcal  A}_{q, 1}$, 
${\mathcal  A}_q$ is a quantum subgroup of $A_{q^m, m}$
for all $m \neq 0$. The quantum groups ${\mathcal  A}_{q, n}$ can be viewed as
a quantum deformation of the quantum groups ${\mathcal  A}_{1, n}$ 
(${\mathcal  A}_{1, 1}$ is just ${\mathcal  A}_1$ 
treated in the previous paragraph.) 
Note that ${\mathcal  A}_q$ is a proper quantum subgroup of
the universal quantum group ${\mathcal  A}$. 

In the degenerate cases, we have 
 $${\mathcal A}_0 = {\mathcal  A}_{0, n}
= {\mathcal  A}_{q, 0} 
= k[a, a^{-1}], \; \; \hfill \text{for all $n \in {\mathbb Z}$},$$
which is commutative and cocommutative. It is clear that 
this is a quantum subgroup of ${\mathcal  A}_{q, n}$ for all $q, n$. 

The underlying algebra of ${\mathcal  A}_q$ above is exactly the same
as the dual algebra $B$ of the
function algebra $A$ of the quantum $E(2)$ group
in Van Daele and Woronowicz \cite{DaeleWor1},
% it sits as a subalgebra of both $A$ and $B$ there
but it has a different Hopf algebra structure from $B$.
The above description is purely algebraic.
It would interesting to see whether one can adopt the $C^*$-version
of $B$ in \cite{DaeleWor1} to obtain a $C^*$-version of ${\mathcal  A}_q$.
For the same reason as in Proposition 3.4 of \cite{W5},
we may not expect the universal
quantum group $\mathcal  A$ to have a $C^*$-version.
\vvvv
%
%{\bf Acknowledgement.}
%Supported by National Science Foundation grant DMS-9627755.
%

%\vspace{1.5cm}
%\hfill %\noindent
%July, 1998

\end{document}